\documentclass[12pt, leqno]{amsart}
\usepackage{polski}
\usepackage{amssymb}
\usepackage{amsmath}
\usepackage[polutonikogreek,   polish, english]{babel}
\usepackage[cp1250]{inputenc}
\date{}
\usepackage{graphicx}
\usepackage{fancyhdr}

\newcommand{\C}{{\mathbb{C}}}
\newcommand{\R}{{\mathbb{R}}}
\newcommand{\Q}{{\mathbb{Q}}}
\newcommand{\N}{{\mathbb{N}}}

\newcommand{\F}{{\mathbb{F}}}

\newcommand{\pv}{\par\vspace{1ex}}
\newcommand{\pvn}{\par\vspace{1ex}\noindent}
\title{A Purely Algebraic Proof of the Fundamental Theorem of Algebra}
\author{Piotr Błaszczyk}
\begin{document}

\fancypagestyle{newstyle}{
\fancyhf{} % clear all header and footer fields
\fancyhead[R]{\thepage}
\fancyfoot[R]{}
\renewcommand{\headrulewidth}{0pt}
\renewcommand{\footrulewidth}{0pt}
}

\pagestyle{newstyle}

\maketitle
\begin{abstract} Proofs of the fundamental theorem of algebra can be divided up into three groups according to  the techniques involved:  proofs that rely on real or complex analysis, {algebraic proofs}, and topological proofs.
 {Algebraic proofs}   make use  of the fact that  odd-degree real polynomials have  real roots.  This assumption, however,   requires analytic methods, namely, the intermediate value theorem for real continuous functions. In this paper, we develop  the idea of algebraic proof further  towards a purely algebraic proof of the intermediate value theorem for real polynomials.   In our proof, we neither  use the notion of continuous function nor  refer to any theorem of real and complex analysis.   Instead, we  apply  techniques of modern algebra: we extend the field of real numbers   to the\linebreak non-Archimedean field of hyperreals {via} an ultraproduct construction and explore some relationships between the  subring of limited hyperreals, its maximal ideal of infinitesimals, and  real numbers.
\end{abstract}

% It is knows an ordered field is real closed if and only if it has an intermediate value property for polynomials.
%
% F is real closed  iff F is euclidean field and every euation of odd degree over F has a root in F. iff F is not algebraically  closed but $\F(\sqrt{-a})$ is algebraically closed.\footnote{See [P.M. Cohn, Algebra, vol. 3, John Wiley \& Sons, Chichester, England 1991], p. 313-315.}
%
% Our proof is  within the framework of real closed field theory.
% we develop further this result and give
% give a \emph{non-analytic}  proof of the intermediate-value theorem for real polynomials -- the one
%  that  To this end we extend the real number field \emph{via} the ultraproduct construction to a non-Archimedean field.  Then, through a combination of techniques usually applied in  finite domains and  ultraproduct, we show that an algebraic equation of odd degree has a real root.
%   Despite this our roof is within the real closed field theory.
%A field is real closed if it is formally real but has no algebraic extension that are formally real

\pvn\textbf {1. Introduction}.
In 1799, Gauss gave the  first widely accepted proof of the fundamental theorem of algebra,  FTA for short: Every nonconstant complex polynomial has a complex root.  Since then, many new proofs have appeared, including new insights, as well as a diversity of tricks, techniques and general methods. Nearly a hundred proofs of FTA  were published up to 1907 (see \cite[p.~98] {RR}).
  Another hundred  proofs
were released in the period 1933 to 2009 (see \cite{BB}). This unusual number of proofs compares  with the multitude of  proofs of the Pythagorean theorem (see \cite{EM}).

   Fine and  Rosenberg \cite{FR}   take a more general,  qualitative perspective and present six exemplary  proofs of  FTA   classified according to the techniques  involved.  They pair these proofs  in accordance with the  basic areas of mathematics and present these pairs  as models of analytic, algebraic and topological proofs.  The first two proofs  require real and complex analysis, the third and  fourth ones apply  algebraic methods: splitting fields and the fundamental theorem of symmetric polynomials,  or the Galois theory and the Sylow theorem. The fifth proof  involves the notion of  the winding number of a closed, continuously differentiable curve $f:\R\rightarrow \C$ around $0$. The sixth one relies on the  algebraic topology  and applies the Brouwer fixed point theorem.\footnote{Some of the recent results
    use linear algebra \cite{HD} and nonstandard analysis \cite{GL}. The first one relies on the intermediate value theorem for real continuous functions, the second one exploits the Brouwer fixed point theorem.}

 After the presentation of analytic and algebraic proofs,  Fine and Rosenberg observe:\\
``We have now seen four different proofs of the Fundamental Theorem of Algebra. The first two were purely analysis, while the second pair involved
a wide range collection of algebraic  ideas. However, we should realize that even in these proofs we did not totally leave \mbox{analysis.} Each of these proofs used the fact the  odd-degree real polynomials have real roots. {This fact is a consequence of the intermediate value theorem, which depends on continuity}. Continuity is a topological property and we now proceed to our final pair  of proofs, which involve \mbox{topology." \cite[p.~134]{FR}}

\pv   All  these proofs mentioned above involve, however,   two kinds of continuity: the continuity of   total order  and the continuity of  function.  The first one is the characteristic feature of real numbers, and since we deal with real polynomials, we cannot ignore the continuity of the reals. The second one refers to  a function -- one can  call the continuity of this kind   a \emph{topological property}.  In our proof of FTA,   the continuity of a function is omitted. In the next section, we develop   this distinction further.

\pvn\textbf {2. Two kinds of continuity.}   The continuity  of the field of real numbers  can be formulated  in many equivalent ways. In this paper, we apply what we believe to be the simplest development -- the one introduced by \mbox{Richard Dedekind}  in  his 1872  \cite{RD}.
To this end, we need a notion of  Dedekind cut.

\pvn\textbf{Definition.} A pair of sets $(L,U)$ is a  Dedekind cut  of a totally ordered set $(X,<)$ if (1) $L,U\neq \emptyset$, (2) $L\cup U=X$,\ (3) $(\forall y\in L)(\forall z\in U)(y< z)$. A cut $(L,U)$ is called a gap if there exists neither a maximum in $L$ nor a minimum in $U$. A cut $(L,U)$ is called a jump    if there exists  both a maximum in $L$ and a minimum in $U$.

\pv  Now,  depending on whether we consider  the real line $(\R,<)$ or  the field of reals $(\R,+,\cdot,0,1,<)$, we get  different continuity axioms.
The categorical characteristics of a continuously  ordered set $(X,<)$, due to George Cantor's \cite{GC}, consists of three conditions: (1) the order is dense, (2) no Dedekind cut of $(X,<)$ is a gap, (3) the order is separable, i.e.,  there exists such a countable subset $Z\subset X$ that is dense in $X$, $$(\forall x,y\in X)(\exists z\in Z )(x<y\Rightarrow x<z<y).\footnote{Next to \cite{GC}, § 11. \emph{Der Ordnungstypus} $\theta$
\emph{des Linearkontinuums X}, see  also \cite{KM}, ch. VI, \S\, 3. \emph{Order types} $\omega, \eta$, \emph{and} $\lambda$.}$$
In other words, any ordered set $(X,<)$ that satisfies the three above conditions is isomorphic to the line of real numbers $(\R,<)$.

The continuity axiom for ordered fields is significantly simpler, for it consists of the requirement (2) alone.  This is because, the density of the field order follows from  ordered field axioms, while  the continuity axiom itself implies  that  the order is  separable (see below).

From now on, $\mathcal F$ denotes a totally ordered field $(\F,+,\cdot,0,1,<)$, that is a commutative field with a total order that is compatible with addition  and multiplication,
\begin{enumerate}
  \item[] $(\forall x,y\in\F)(x<y\Rightarrow x+z<y+z)$,
  \item [] $(\forall x,y,z\in \F)(x<y, 0<z \Rightarrow x\cdot z<y\cdot z)$.
\end{enumerate}

 %accordingly, by $\mathcal R$ we mean the real number field \mbox{$(\R,+,\cdot,0,1,<)$}.

\pvn\textbf{Definition.} The field of real numbers   is  an  ordered field $\mathcal F$, in which
 every Dedekind cut $(L,U)$ of $(\F,<)$ satisfies the condition
\begin{equation}\tag{C1}(\exists x\in\F)(\forall y\in L)(\forall z\in U)(y\leq x\leq z).
\end{equation}

The  categoricity theorem  states that any two ordered fields satisfying axiom (C1) are isomorphic (\cite[p. 105]{CE}). In other words,
any ordered field  satisfying (C1) is isomorphic to the field of real numbers $(\R,+,\cdot,0,1,<)$.

The well-known constructions of the reals, e.g., the one that identifies real numbers with  cuts of the line of rational numbers $(\Q,<)$  due to  Dedekind, show that there exists at least one field of real numbers. On the other hand,  the  categoricity theorem implies that there exists at most one, up to isomorphism, field of real numbers.

 The order of  an ordered field    is dense: it follows from the simple observation that if $x<y$, then $(x+y)/2$ lies between $x$ and $y$. On the other hand, the density of \mbox{a field} order is equivalent to the claim that no cut of $(\F,<)$ is a jump. As a result, (C1) is reduced to  the fact that no Dedekind cut of $(\F,<)$ is a gap.\footnote{It follows from the axiom (C1) that the set of fractions, $\Q$, is dense in $(\R, <)$; in fact, it is an equivalent form of the Archimedean axiom (see \S\, 4 below). Thus, for ordered fields, (C1) suffice to characterize the continuity of the order.}  With this knowledge we can easily visualize the continuity of the real line.
 \newpage
\begin{center}\includegraphics[scale=0.4]{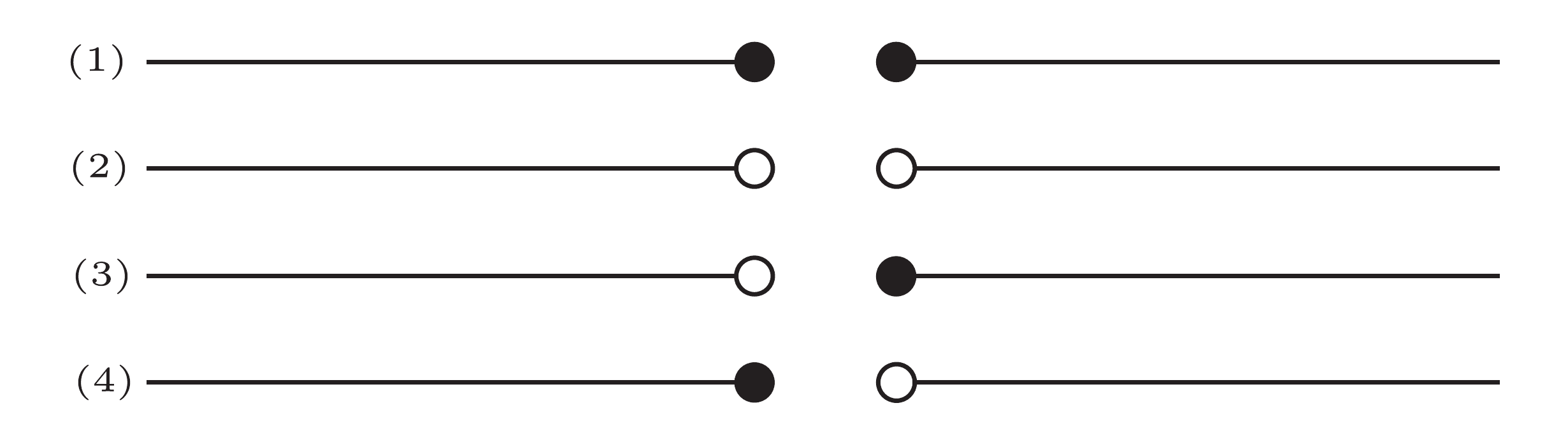}

{\tiny\ \   Dedekind cuts}\end{center}

On the above diagram a black dot corresponds to the maximal element of the class $L$, or the minimal element of the class $U$. If  there are,  accordingly,  no maximum or minimum, then we put a  white dot.   In this way the diagram \emph{Dedekind cuts} represents the only possible kinds of cuts of the line $(\F,<)$.  Cut (1) is a jump,  cut (2) is a gap.  Axiom (C1) states that  cuts (3) and (4) are the only possible cuts on the line of real numbers $(\R,<)$.

  (C1) is called   the continuity axiom. In fact, it is one of the many equivalent formulations of  the continuity of the reals. Here are some other popular versions:
\pvn (C2) If $A \subset \F$ is a nonempty set which is bounded above, then there exists $a \in \F$ such  that  $a = \sup A$. \\
(C3) Field $\mathcal F$ is Archimedean and every Cauchy (fundamental) sequence $(a_n)\subset \F$  has a limit in $\F$.

\pv Cohn and Ehrilch \cite[p. 95--96]{CE},  present other standard formulations  of the continuity axiom  which appeared at the turn of 19th and 20th century.  Due to their existential character (``there exists $x\in\F$ such that"), they all  make it possible to  determine specific real numbers -- each version with a different technique, e.g.,  via  Dedekind cuts,  bounded sets, or fundamental sequences.
In a sense, each form of the continuity axiom corresponds to a technique. In  the next section we  present  a new version of the continuity axiom that  reflects the technique we adopt  in our proof of FTA.

\pvn \textbf{2.1.} Now we  turn to the  continuity of function.
\pvn \textbf{Definition.} A function $f:\F\mapsto \F$ is continuous at a point $a\in \F$  if the following
 condition holds

$(\forall \varepsilon>0)(\exists \delta>0)(\forall x\in\F)(|x-a|<\delta\Rightarrow |f(x)-f(a)|<\varepsilon).$

\pv To elaborate, first, observe that the definition of a continuous function can be formulated  in any ordered field.   Similarly, we can develop the theory of limits of sequences in any order field  (see \cite{CE}, chap. 3). Moreover, in  any ordered field we can also  formulate the intermediate value theorem (IVT for short).

Secondly,   the standard proof of  IVT  for a real continuous function proceeds as follows:
Let  $f:\R\mapsto\R$ be a continuous function with $f(0)<0$ and $f(1)>0$. Putting $A=\{x\in[0,1]: f(x)<0\}$,
by (C2), we can take $a=\sup A$. The continuity of  $f$  implies $f(a)=0$.  The final step in this argument is based on the rule stating  that a real continuous function   preserves the sign:
\begin{equation}\tag{SR}f(c)\neq 0\Rightarrow (\exists \delta>0)(\forall x\in \R)(|c-x|<\delta\Rightarrow f(c)\cdot f(x)>0).\end{equation}

Thus, if $f(a)< 0$, then, by (SR), for some $\delta>0$ and every $x$  such that  $0<x-a<\delta$ holds $f(x)<0$, contrary to the assumption $a=\sup A$; the same argument applies to  the case  $f(a)>0$. Hence, via the trichotomy law of the total order,  $f(a)$ equals $0$.

This proof clearly manifests the combination of the continuity of the reals and the continuity of a function: we apply the (SR) rule to the point defined by (C2).\footnote{Essentially, this argument  goes back as far as  Bolzano \cite{GB}.}

Thirdly, from the logical point of view, the difference between the continuity of an  ordered field  and the continuity of  function is not an absolute one. Axiom (C1) turns out to be equivalent to many  statements of real analysis which are typically presented as theorems and involve the notion of a function. Teismann  \cite{HT} proves that IVT, as well as  the mean value theorem,  is equivalent to (C1).\footnote{The proof of the first equivalence  can  be easily  based on our diagram \emph{Dedekind cuts}. If a field $\mathcal F$ does not satisfy (C1), then there is a cut of $(\F,<)$ which yields a gap $(L,U)$, as the one marked by (2). A function given by $f(x)=0$, for $x\in L$, and $f(x)=1$, for $x\in U$, is continuous and does not satisfy IVT. If $\mathcal F$ is the real number field, then we can  adopt  the proof presented above to show  that a real continuous function satisfies IVT.}
Riemenschneider \cite{OR}  lists 37 versions of the continuity axiom, mostly theorems of one-variable real analysis.

\pvn \textbf{3. Real closed fields.} The theory of real closed fields, started by Emil Artin and Otto Schreier in the 1920s, provides a general framework for this paper.  Our sketchy account of this theory  starts with a technical notion of a formally real field.\footnote{The general reference here is \cite{AS},  \cite{PC}, chap. 8, \cite{NJ}, chap. 6;  see also  \cite{HZ} as it applies the theory of real closed field to give a constructive proof of FTA.}
At the end of this section we show that in order to prove  FTA it  suffices to prove IVT for polynomials.
\pvn\textbf{Definition.} A field $(\F,+,\cdot,0,1)$ is formally real if the sum of squares is zero only if  each summand is zero.\footnote{We assume that all fields considered are commutative.}

\pv An equivalent formulation of this definition is this: A field is formally real if  $-1$ cannot be written as a sum of squares, that is
\mbox{$-1\notin\big\{\sum a_i^2:\ a_i\in\F\big\}$.}
Thus, any ordered field is  obviously formally real.

\pvn\textbf{Definition.} A field $(\F,+,\cdot,0,1)$ is real closed if it is formally real and  every proper algebraic extension of the field is not formally real.

\pv A real closed field is in fact an ordered field. Setting
\[x<y \Leftrightarrow y-x\in \big\{\sum a_i^2:\ a_i\in\F\big\}\]
we obtain a total order  compatible with the addition and multiplication \mbox{on $\F$.}

 The real number field and the real algebraic number field are both really closed, so the order of a real closed field is not necessarily  continuous.
Next, some theorems  concerning real polynomials, like the intermediate value theorem,   mean value theorem,  extreme value theorem,  Rolle's theorem, or  theorem of Sturm about the number of zeros in an interval,  also hold in a real closed field.\footnote{See \cite{NJ}, chap. VI, \S\, 3.} These results suggest  that the algebraic features of the real number field alone \mbox{imply IVT.}

\pvn \textbf{3.1.} The theory of real closed fields provides a criterion  for a field to be algebraically closed. Artin and Schreier \cite{AS} show that: ``In a real closed field, every polynomial of odd degree has at least
one  root."\cite[p.~ 275]{AS}

Then, they proceed to prove the  proposition that constitutes the basis for our argument, namely:
\pvn    ``A real closed field is not algebraically closed. On the other hand, the field obtained by adjoining $i\ (=\sqrt {-1})$   is algebraically closed." \cite[p. 275]{AS}\footnote{The converse of this  theorem also holds:
 If $\mathcal F$ is such an ordered field  that the field obtained by adjoining $\sqrt{-1}$ is algebraically closed, than $\mathcal F$ is real closed (see \cite[p. 277]{NJ}).}

    \pv In fact,  Artin and Schreier present two  proofs which make the essence of what  Fine and Rosenberg call the algebraic proofs of FTA.

 The above result supports the conclusion that  FTA is equivalent to the claim  that real numbers form a real closed field.

\pvn \textbf{3.2.} The definition of a real closed field does not provide  simple  criteria for a field to be really closed.
Yet, real closed fields can be described more explicitly.
Cohn \cite{PC} shows that an ordered field $(\F,+,\cdot,0,1,<)$ is real closed if and only if it is closed under a square root operation
and   odd-degree polynomials  have  roots in $\F$, that is
\begin{enumerate}
\item[(R1)] $(\forall x>0)(\exists y\in \F)(y^2=x)$,
 \item [(R2)] $(\forall a_1\in \F)... (\forall a_{2n}\in \F)(\exists x\in \F)(x^{2n+1}+\sum_{i=1}^{2n}a_ix^i=0)$.
 %If $f\in\F[x]$,  $deg f =2n+1$,  then $(\exists y\in \F)(f(y)=0)$.
\end{enumerate}

Indeed, these properties, along with  ordered field axioms, constitute the set of  \mbox{axioms} for ordered real closed fields (see \cite[p. 94--95]{DM}).\footnote{Given $(\F,+,\cdot,0,1)$ is formally real, the axiom (R1) has to take the following form: $(\forall x\in \F)(\exists y\in \F)(y^2=x\vee y^2=-x)$.}
It is also easily seen  that  properties (R1), (R2) can be combined into  this one:
\pvn ``An ordered field is real closed if and only  if it  has the intermediate value property  for polynomials." \cite[p. 315]{PC}

\pv As a result we obtain:

\pvn\textbf{Proposition 1.} If an ordered field  has  the intermediate value property for polynomials, than the field obtained by adjoining $\sqrt {-1}$   is algebraically closed.

\pv Thus, to  prove  FTA it is sufficient to prove IVT for  real polynomials.

\pvn\textbf{4. Archimedean fields.} In the next section we will deal with a non-Archimedean field,   so now  we restate the definition  and some basic facts concerning  Archimedean   fields.

\pvn \textbf{Definition.} A totally ordered field $\mathcal F$ is Archimedean if it
satisfies the condition
\begin{equation}\tag{A1}(\forall x\in \F)(\exists n\in\N)(n>x).
\end{equation}

Axiom (A1) is called the Archimedean axiom. Here is its equivalent version
\begin{equation}\tag{A2}\lim\limits_{n\to \infty} \frac 1n=0.
\end{equation}

 Yet another equivalent formulation of the Archimedean  axiom  is the following one:
 If  $(L,U)$ is a Dedekind cut of $(\mathbb{F},<)$, then
\begin{equation}\tag{A3} (\forall n\in\mathbb{N})(\exists x\in L)( \exists y\in U) (y-x< 1/n).\end{equation}

 Given $(L,U)$ is a Dedekind cut of $(\F,<)$, by (A3), we  can find  such a sequence $(r_n)\subset\F$  that satisfies  conditions
\begin{equation} r_{2k-1}\in{L},\ r_{2k}\in{U},\ r_{2k}-r_{2k-1}<\frac{1}{k},\end{equation}
\[r_1\leq r_3\leq...\leq r_{2k-1}\leq ...\leq r_{2k}\leq ...\leq r_4\leq r_2.\]

We  refer to this fact proving proposition 3 below.

The real number field is the \emph{biggest} Archimedean field, that is any Archimedean field can be embedded into the field of reals; on the other hand,  any extension of the reals  is a non-Archimedean field (see \cite{CE}).

\pvn \textbf{5. Extending the  real number field}.
In this section we  introduce the construction called the ultraproduct  to extend the field of reals to the field of hyperreals (nonstandard real numbers). We apply this construction directly to the field of reals to build a special kind of ultraproduct called an ultrapower. This presentation follows \cite{PC}, \cite{RG}.

To start with,  we give the definition of an ultrafilter.
  \pvn \textbf{Definition.} A family of sets $\mathcal U\subset \mathcal P (\N)$ is an ultrafilter on $\N$ if (1) $\emptyset \notin \mathcal U$, (2) if $A,B\in\mathcal U$, then $A\cap B\in\mathcal U$, (3) if $A\in\mathcal U$ and $A\subset B$, then $B\in\mathcal U$, (4) for each  $A\subset \N$, either $A$ or its complement $\N\setminus A$ belongs to $\mathcal U$.

It follows from this definition that either the set of odd numbers or the set of even numbers belongs to an ultrafilter.

  Now, take the family of  sets with  finite complements. It obviously  satisfies conditions (1)--(3) listed in the definition of an ultrafilter. By  Zorn's lemma, this family  can be  extended to an  ultrafilter $\mathcal U$  on $\N$ (see \cite[p. 29]{PC}, \cite[p. 20--21]{RG}).

  From now on, $\mathcal U$ denotes a fixed ultrafilter  on $\N$ containing  every  subset with a finite complement.\footnote{The reader can find the same reasoning in the literature  encoded in  one short sentence: Let $\mathcal U$ be  a fixed nonprincipal ultrafiltr on $\N$.}

 In the product $\R^{\N}$ we define a relation
\[(r_n)\equiv (s_n)\Leftrightarrow \{n\in \N:\ r_n=s_n\}\in\mathcal U.\]
This is easily seen to be an equivalence relation, satisfying reflexivity, symmetry and transitivity.
Let $\R^*$  denote  the reduced product $\R^{\N}/_\equiv$.

The equality relation in $\R^*$  is obviously given by
\[[(r_n)]= [(s_n)]\Leftrightarrow \{n\in \N:\ r_n=s_n\}\in\mathcal U.\]

It follows from the notion of ultrafilter that
\[[(r_n)]\neq [(s_n)]\Leftrightarrow \{n\in \N:\ r_n\neq s_n\}\in\mathcal U.\]

 Algebraic operations on $\R^*$ are defined pointwise, that is
\[[(r_n)]+[(s_n)]=[(r_n+s_n)],\ \ \  [(r_n)]\cdot [(s_n)]=[(r_n\cdot s_n)].\]

And a total order on $\R^*$ is given by the following definition
\[[(r_n)]<[(s_n)]\Leftrightarrow \{n\in \N :r_n<s_n\}\in\mathcal U.\footnote{In the above definitions, we adopt a standard convention to use the same signs for the relations and operations  on $\R$ and $\R^*$; see \cite{RG}, chap. 3,\linebreak \S\, 3.6.  }\]

We embed the set of reals $\R$  into the set of hyperreals $\R^*$ by  identifying the standard real number $r$ with the hyperreal determined by the constant sequence $(r,r,r,...)$. Thus,  setting   $r^*=[(r,r,r,...)]$, the embedding is given by  the following map
\[\R\ni r\mapsto r^*\in \R^*.\]

To end this part of our development,  we would like to ease our notation: from now on  we  use $r$ for    hyperreal  number $r^*$.
In fact, it is a consequence of the  convention that  $\R$ is a \emph{subset} of $\R^*$.

\pv\textbf{Proposition 2.} $(\R^*,+,\cdot,0,1,<)$ is non-Archimedean ordered field.

{\sc Proof}. Goldblatt \cite[p. 23--24]{RG}  gives a straightforward proof that the hyperreals form  a totaly ordered field.\footnote{This claim is also a straightforward consequence of  the Łoś theorem \cite{JL},  also known as  the transfer principle; however,  in this paper we avoid arguments relying on mathematical logic.}
In addition, %\linebreak
  by (A2), for any $r\in \R$ the inequality holds $r\neq [(1/n)]$, which means that the field of  hyperreals   extends  the field of reals. Thus, the  field of hyperreals  is   non-Archimedean field.
\hfill{$\Box$}

\pv In the last section,  we will  also  show that the field of hyperreals  is  real  closed.
\pvn \textbf{5.1.} We define, in a standard way,  subsets of $\R^*$ --  sets  of infinitesimals,   limited and infinitely large hyperreals, namely\footnote{Our definitions agree with the those given by Artin and Schreier \cite{AS}. However, we can offer some simple examples of infinitesimals and infinitely large numbers, namely $[(1/n)]\in \Omega$, $[(n)]\in\Psi$, moreover, it is easy to demonstrate that if \mbox{$\lim\limits_{n\rightarrow\infty}r_n=0$,} then $[(r_n)]\in\Omega$, and if $\lim\limits_{n\rightarrow\infty}r_n=\infty$, then $[(r_n)]\in\Psi$.}
\begin{eqnarray*}
x\in \Omega      &\Leftrightarrow& (\forall n\in \N)(|x|<1/n),\\
x\in \mathbb L &\Leftrightarrow & (\exists n\in \N)(|x|<n),\\
x\in \Psi             &\Leftrightarrow & (\forall n\in \N)(|x|>n).\end{eqnarray*}

 The set of  positive  hyperintegers (hypernaturals) $\N^*$ is defined by
 $$\N^*=\{[(r_n)]\in\R^*:\  \{n\in\N:\ r_n\in\N\}\in\mathcal U\}.$$

 Roughly speaking, the set $\N^*$ consists of elements $[(n_j)]$, where $(n_j)\subset\N$.

Next,  on the set $\R^*$ we define  a  relation    \emph{$x$ is infinitely close to $y$} by putting
\[x\approx y \Leftrightarrow x-y\in\Omega.\]

This is easily seen to be an equivalence relation, satisfying reflexivity, symmetry and transitivity.

  Following are some elementary facts concerning  these concepts.  For the sake of completeness, we present  a short justification for each one, though they are almost obvious.

\pvn (F1) A standard real is a limited hyperreal,  $\R\subset \mathbb L$.

It is the consequence of the Archimedean axiom.

\pvn (F2) Two standard real numbers $r, s$ do not lie infinitely close to each other.

It follows from the Archimedean axiom that  the real  number
$|r-s|$ is greater than $1/k$ for some $k\in \N$. Thus this number is not infinitesimal, and  neither $r-s$ nor $s-r$ belongs to $\Omega$.

\pvn (F3) Limited numbers form  an ordered ring $(\mathbb L,+,\cdot,0,1,<)$  with $\Omega$ being its maximal ideal. Particulary, the following condition is satisfied
\[(\forall x \in \Omega)(\forall y\in\mathbb L)(x\cdot y\in\Omega).\]

The first part of this claim  is a consequence of the interplay  between the quantifiers ``for all"  and ``exists" occurring in the definitions of  sets $\Omega$ and $\mathbb L$. To show that $\Omega$ is the maximal ideal of the ring $\mathbb L$, suppose, to obtain a contradiction, that
$G$ is an ideal of $\mathbb L$ such that $\Omega\varsubsetneq G$ and $G\varsubsetneq \mathbb L$. Take $x\in G\setminus \Omega$. Since $x$ is limited, for some $m\in\N$   holds $|x|<m$; since it is not infinitesimal, for some $k$ holds $|x|>1/k$.   Then
\[  1/k<|x|<m.\]
Hence, via the rules of an ordered field, we obtain
$$1/m<|x^{-1}|<k,$$ which means that $x^{-1}$ is a limited hyperreal. Since $x$ belongs to the ideal $G$, the element $1=x\cdot x^{-1}$ also belongs to $G$, a contradiction.\footnote{By proposition 3, one can show that the quotient ring $\mathbb L/\Omega$ is isomorphic to the real number field. As a result, we can represent the set  $\mathbb L$  as the sum of disjoint sets $r+\Omega$, the so-called \emph{monads}, for $r\in \R$.}
\pvn (F4) If $x\in \R^*$ and $x\neq 0$, then the equivalence holds
\[x\in \Omega \Leftrightarrow x^{-1}\in\Psi.\]

If
 $x\in \Omega\setminus\{0\}$  and  $x^{-1}\notin \Psi$, then $x^{-1}\in \mathbb L$.  By (F3),  $x\cdot x^{-1}\in \Omega$, contrary to the fact that $1$ is not an infinitesimal. Next, if $|x^{-1}|>n$, for every $n\in \N$, then $x<1/n$, for every $n\in\N$. It is equivalent to the claim that if $x^{-1}$ is infinitely large, then $x$ is infinitely small.

To summarize this subsection, we present a diagram representing the ultraprower construction.

\begin{center}\includegraphics[scale=0.6]{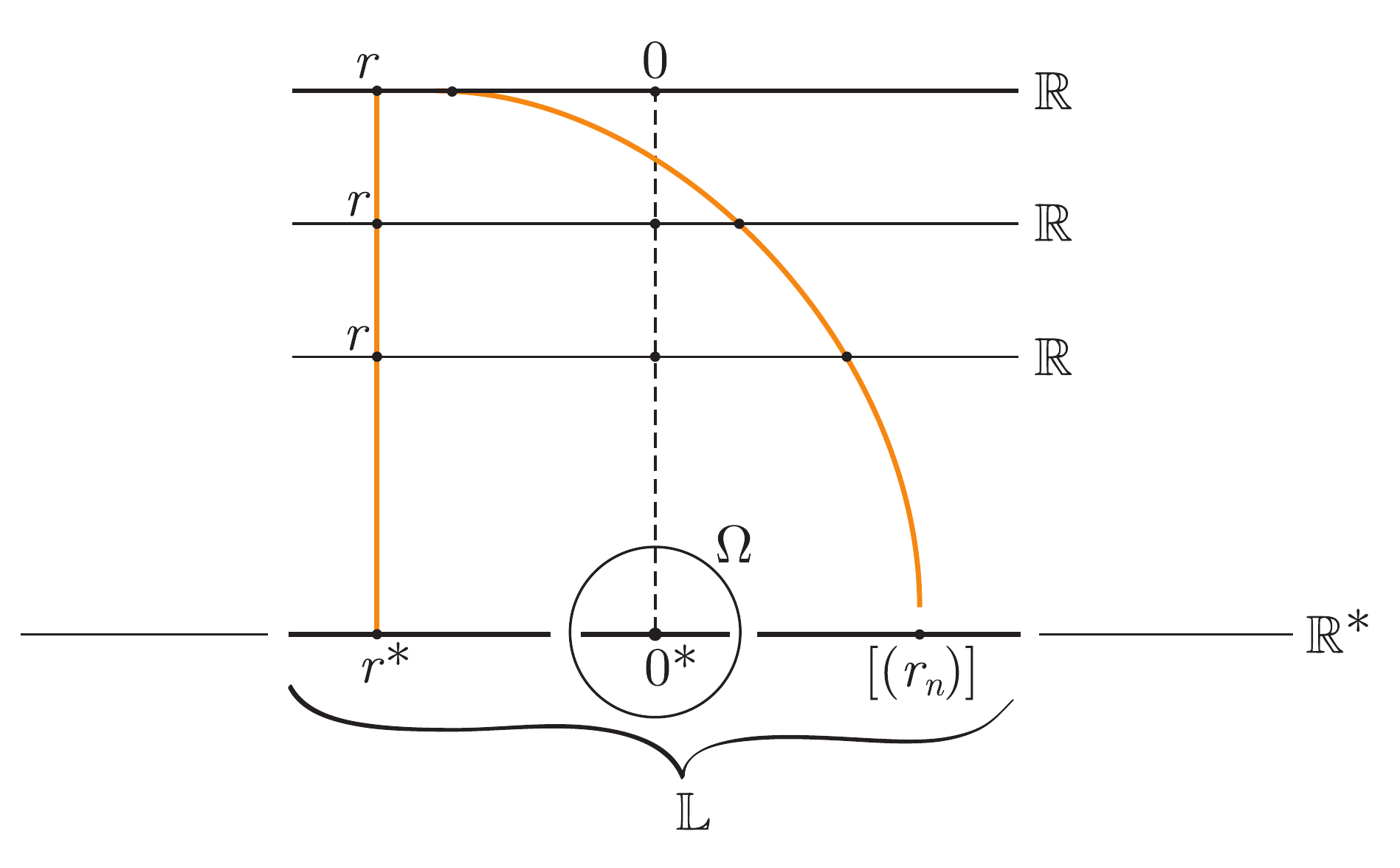}\end{center}

\pvn \textbf{5.2.}   We can apply the ultraprower construction to any ordered \mbox{field $\mathcal F$.} Thereby we obtain the set
$\F^*$ and its subsets $\Omega$, $\mathbb L$, $\Psi$, as well as the relation $x$ \emph{is infinitely close to} $y$.  With these notions
we can formulate yet another version of the continuity axiom, namely
\begin{equation}\tag{C4}(\forall a\in\mathbb L)(\exists! z\in\mathbb{F}) (a\approx z).
\end{equation}

\pvn {\bf Proposition 3.}  The statements (C1) and (C4) are equivalent.

\pvn {\sc Proof}.  The first part of this claim -- if $\mathcal F=(\R,+,\cdot,0,1,<)$,  then  each limited number $a$ is infinitely close to exactly one real number $z$ --  is the well-known Standard Part Principle (see  \cite{JB},  \cite[p. 53]{RG}). In the  proof that follows, we  apply the Archimedean axiom in the version  (A3).

 The limited number  $a\in\R^*$  determines a Dedekind cut of $(\R,<)$
\[L=\{x\in \R: x\leq a\},\ \ \ U=\{x\in \R: x>a\}.\]

By (C1), the cut $(L,U)$ determines the real number $z$. By (1), we find a sequence $(r_n)\subset\R$ such that
 $r_{2k-1}\in L$, $r_{2k}\in U$ and $r_{2k}-r_{2k-1}< 1/k$. Hence,
\[z,a\in \bigcap\limits_{k=1} [r_{2k-1},r_{2k}], \]
which gives  $z\approx a$.

  Since two standard reals do not  lie infinitely close to each other,   only one   real    is infinitely close to $a$. This unique number is called the standard part, or the {shadow},  of $a$ and is denoted by $a^o$. Thus $a^o\approx a$, or $a=a^o+\varepsilon$, for some $\varepsilon\in\Omega$.
Also, note that if $a\approx b$, then, by the uniqueness of the standard part, we obtain $a^o=b^o$.

%To show our claim suppose that $[(q_n)]$ is not infinitesimally close to $a=[(a_n)]$. Then,
%$\{n\in \N |\ |a_n-q_n|\geq 1/k \}\in{\mathcal{U}}$, for some $k\in \N$. Moreover,
%$[(q_n)]\prec a$ or $[(q_n)]\succ a$, either $\{n\in \N|\ a_n\in{A}\}\in \mathcal{U}$ or $\{n\in \N|\ a_n\in{B}\}\in\mathcal{U}$.
%\begin{equation}\tag{C2}(\forall x\in \mathbb L)(\exists! y\in \F)(x\approx  y^*).\end{equation}

For the second part of the proof, we  first show that (C4) implies (A2).\footnote{For  another  proof of  proposition 3  see \cite{Ha}, however, it implicitly relies on the assumption that $\mathcal F$ is Archimedean.}     Suppose, on the contrary,  that (A2) does not hold.
Then, for some $\varepsilon>0$ and for every $n$ holds $\varepsilon<1/n$. Hence  $0\approx [(1/n)]\approx \varepsilon$, which contradicts the claim  that there is only one element in $\F$ infinitely close to  $[(1/n)]$.

Thus  we come to the main part of the proof.     Let $(L,U)$ be a Dedekind cut of $(\F,<)$
and  $(r_n)\subset \F$  satisfy condition (1). Then
\mbox{$r_1<[(r_n)]<r_1+1$}. Since $\mathcal F$ is Archimedean, for some $n\in\N$ holds $[(r_n)]<n$. As a result we obtain that $[(r_n)]\in \mathbb L$.  By (C4), there exists $z\in \F$  such that $z\approx [(r_n)]$. We show that $z$ is the greatest element in $L$, or the least element in $U$.

Seeking a  contradiction, suppose that $(L,U)$ is a gap.
 We need to consider four possibilities resulting from a combination of the following  conditions: (1)  \mbox{$z\in{L}$,} (2) $z\in U$, (a) the set of odd  numbers belongs to $\mathcal U$, (b) the set of even numbers belongs to $\mathcal U$.

(Ad 1a.)  Suppose  $z\in{L}$ and  the set of odd numbers belongs to $\mathcal U$.  There exists $x\in{L}$, such that $z<x$, for  $(L,U)$ is a gap.
Set
 \begin{equation}\theta=\frac{x-z}{2}.\end{equation}

Since   $z\approx [(r_n)]$, it follows that
 \[\left\{n\in\mathbb{N}\ |\ |r_{n}-z|<\theta\right\}\in\mathcal{U}.\]

Put
\[ A=\{n\in\N: n\ \mbox{is\ odd}\}\cap\{n\in\mathbb{N}\ |\ |r_{n}-z|<\theta\}\cap\{n\in\N:\frac{1}{n}<\theta\}.\]

First, $A\in \mathcal U$. Second, if $k\in A$, then $r_{k}-z<\theta$; in consequence
\begin{equation}r_{k}<x-\theta.\end{equation}

By (2), we have $r_k<x$. Thus $r_{k}\in L$.

 On the other hand, it follows from (1) that
  \begin{equation}r_{k+1}-r_{k}<\frac{1}{k}<\theta.\end{equation}

 By adding inequalities (3) and (4)
we obtain $r_{k+1}<x$. Hence $r_{k+1}\in{L}$. But $k+1$ is even, so $r_{k+1}\in U$,  contrary to the assumption $L\cap U=\emptyset$.

(Ad 1b.) In the same manner, suppose  $z\in{L}$ and  the set of even numbers belongs to $\mathcal U$.  Let $x\in L$ be such that $z<x$; set
  $\theta=\frac{x-z}{2}$.
 Put
\[ A=\{n\in\N: n\ \mbox{is\ even}\}\cap\{n\in\mathbb{N}\ |\ |r_{n}-z|<\theta\}\cap\{n\in\N:\frac{1}{n}<\theta\}.\]

  Let $k\in A$.  As before,
  $r_k<x$, which gives  $r_k\in L$. Since $k$ is even, the term $r_k$ belongs to $U$. Thus $r_k\in L\cap U$, a contradiction.

  The same reasoning applies to  cases (2a), (2b).  \hfill{$\Box$}

\pvn \textbf{5.3.} Let $r\in \R$. In our proof of FTA, we will  also need these simple facts:\footnote{One  can  consider these facts  a nonstandard counterpart of the rule (SR),  given above in section 2.}

If $r>0$, then for  every $\varepsilon\in\Omega$ the relation obtains $r+\varepsilon>0$.

 If  $r<0$, then for every $\varepsilon\in\Omega$ the relation obtains $r+\varepsilon<0$.

For the proof of the first fact, note that if $r>0$, then for $\varepsilon \in \Omega$ holds $|\varepsilon|<r/2$. Hence $r+ \varepsilon >r/2$, and $r+\varepsilon>0$.

These two facts imply the third result:

If $x\approx y$ and $x\cdot y<0$, then $x,y\in\mathbb L$ and the standard part of $x$ is equal to $0$,
\begin{equation}\label{rz}(\forall x,y\in\mathbb R^*)(x\approx y, x\cdot y<0\Rightarrow x^o=0). \end{equation}
\pvn \textbf{6. Polynomials.} Let $f\in \R[x]$ be a real polynomial,
\[f(x)=a_0+a_1x+...+a_mx^m,\ \  \mbox{where}\ \ a_i\in \R.\]

By  $f^*$ we mean a hyperreal polynomial
  $f^*:\R^*\mapsto \R^*$ with real coefficients $a_i$,
\[f^*(x)=a_0+a_1x+...+a_mx^m,\]
 defined by
\begin{equation}\label{ff}f^*([(r_n)]) =[(f(r_1),f(r_2),...)].\end{equation}

If $r\in \R$, then $f^*(r)=[(f(r),f(r),...)]$.
  Since we identify real number $r$ with hyperreal $r^*$,  the equality $f^*(r)=f(r)$ obtains.

\pvn \textbf{Lemma.} Let $f\in\R[x]$ and $a\in\mathbb L$. If $a\approx b$, then
 $f^*(a)\approx f^*(b)$.

{\sc Proof.} Set
\[f^*(x)=a_0+a_1x+...+a_mx^m, \ \ \  a_i\in \R.\]

Let $r$ be the standard part of $a$, that is $r=a^o$. Thus $r\approx a$, and for some $\varepsilon\in \Omega$ we have   $a=r+\varepsilon$.

Now, for the real number $r$, the following equalities hold
\begin{eqnarray*}f^*(r+\varepsilon)&=&a_0+a_1(r+\varepsilon)...+a_m(r+\varepsilon)^m\\
                                                         &=& f^*(r)+ \varepsilon\cdot  w(a_1,...,a_m,r,\varepsilon),
                                                         \end{eqnarray*}
   where $w(a_1,...,a_m,r,\varepsilon)\in \mathbb L$.  Since infinitesimals form an ideal of the ring $\mathbb L$,  the hyperreal number $\varepsilon\cdot  w(a_1,...,a_m,r,\varepsilon)$ belongs to $\Omega$.   Hence
\mbox{$f^*(r+\varepsilon)-f^*(r)\in\Omega$.} The identification $f^*(r)=f(r)$ clearly forces $f^*(r+\varepsilon)\approx f(r)$. Since $a=r+\varepsilon$, we obtain
\[f^*(a)\approx f(r).\]

If $a\approx b$, then via the transitivity of the relation  \emph{is infinitely close}, we get $r\approx b$.
By the uniqueness   of the standard part,  we also have $r=b^o$.
The  reasoning applied to the pair $r,a$ works for the pair $r, b$ too, thus
\[f^*(b)\approx f(r).\]

Finally,  once more applying the transitivity of the relation  \emph{is infinitely close}, we have
\[f^*(a)\approx f^*(b).\]
 \hfill{$\Box$}

One can consider this lemma  as a nonstandard counterpart of the standard claim that a real polynomial is a continuous function. Indeed, Birkhoff and Mac Lane, in Chapter 4, entitled \emph{Real Numbers}, of their \cite{BMc}  provide a proof for this claim which applies the  $\varepsilon-\delta$ technique. In the next chapter, introducing   the proof of  FTA  they write: ``Many proofs of this celebrated theorem are known. All proofs involve nonalgebraic concepts like those introduced in Chap 4" \cite[p. 114]{BMc}.

\pv We  consider  the proof of the following proposition       purely algebraic. By Proposition 1, this proposition  is equivalent to  FTA.
\pvn \textbf{Proposition 4.} Let $f\in \R[x]$ and  $a,b\in\R$. If $f(a)\cdot f(b)<0$, then for some $c\in(a,b)$ holds $f(c)=0$.

\pv {\sc Proof.}  Obviously, we can   take $a=0$, $b=1$.
Suppose $f(0)<0$ and $f(1)>0$.
By (\ref{ff}) the  same relation obtains in the realm of hyperreals, that is  $f^*(0)<0$ and $ f^*(1)>0$.

Let  $N=[(n_j)]=[((n_1, n_2,...)]$ be an infinitely large hyperinteger; we can take, for example,  $N=[(n)]$.  Set
\[I_N=\{K/N:  0\leq K\leq N\}=\Big\{0,\frac{1}{N}, \frac{2}{N},...,\frac{N-1}{N}, 1\Big\}.\]

In a similar way we  define  sets $I_{n_j}$,
\[I_{n_j}=\{k/n_j:  0\leq k\leq n_j\}=\Big\{0,\frac{1}{n_j}, \frac{2}{n_j},...,\frac{n_j-1}{n_j}, 1\Big\}.\]

However, while sets $I_{n_j}$ are finite, the set $I_N$ is infinite; in fact, it has cardinality continuum.\footnote{The set $I_N$ is usually called a  \emph{hyperfine grid}.}

The image  of $I_{n_j}$ under the map $f$  is a finite set
\[f(I_{n_j})=\{f(0),f(1/n_j),..., f(1)\}.\]

If for some $k$ we have $f(k/n_j)=0$, then the proof is  done.  Thus, we can assume  the elements  of  $f(I_{n_j})$ are either negative or positive.
Let $k_j$, where $0\leq k_j< n_j$, be the first integer such that  $f(k/n_j)<0$ and  \mbox{$f((k+1)/n_j)>0$},
\begin{equation}\label{fin} f(k_j/n_j)<0,\ f((k_j+1)/n_j)>0 .\end{equation}

We can always find such $k_j$, for  the sets $I_{n_j}$ are finite.

Put $K=[(k_j)]$.  Since $0\leq k_j<n_j$, hyperintiger $K$ satisfies the inequalities $0\leq K<N$.  Moreover, $K+1=[(k_j+1)]$.

By (\ref{ff}) and (\ref{fin}) we have
\begin{equation} f^*(K/N)<0,\ \ f^*((K+1)/N)>0. \end{equation}

Number $K/N$ is  a limited hyperreal. Let $c$ be its standard part, that is $c=(K/N)^o$.

   Since $N$ is infinitely large, the element  $1/N$ belongs to $\Omega$.
  Thus $K/N\approx (K+1)/N$, and via the transitivity of  the relation \emph{is infinitely close}  we obtain
  \begin{equation} K/N\approx c\approx (K+1)/N. \end{equation}

Next,  by the lemma, it follows  that
\begin{equation}  f^*(K/N)\approx f(c)\approx f^*((K+1)/N).\end{equation}

Now we come to the final part of the proof. By  (10) hyperreals $f^*(K/N)$ and $f^*((K+1)/N)$ are  infinitely close;  by (8)  they have opposite signs. Then,  by (\ref{rz}),
 the standard part of $f^*(K/N)$ is equal to $0$, that is
 \[(f^*(K/N))^o=0.\]

 On the other hand, by Proposition 3, there exists one and only one standard real number infinitely close to $f^*(K/N)$; by (10)  it is the number $f(c)$, that is
 \[(f^*(K/N))^o=f(c).\]

  Hence
then $f(c)=0$.

To end the proof, observe that the hyperreal $K/N$ belongs to the segment $(0,1)$, so its standard part, $c$, lies in the segment $[0,1]$. Since
\mbox{$f(0)\cdot f(1)<0$,} the real number $c$  equals neither $0$ nor $1$. As a result, $c$  lies in the segment $(0,1)$.
\hfill{$\Box$}

\pv Following, we show that a  polynomial with hyperreal coefficients  has the intermediate value property.
\newpage
\pvn \textbf{Proposition 5.} Let $f\in \R^*[x]$ and  $a,b\in\R^*$. If $f(a)\cdot f(b)<0$, then for some $c\in(a,b)$ holds $f(c)=0$.

{\sc Proof.} Let a hyperreal polynomial
\[f(x)=A_0+A_1x+...+A_mx^m,\ \ \ A_i\in\R^*,\ m\in\N,\]
be such that $f(a)<0$ and $f(b)>0$, with  $a=[(a_n)]$ and $b=[(b_n)]$.
Suppose   $A_i=[(r_{i,n})]$, where $0\leq i\leq m$. The hyperreal polynomial  $f$  is accompanied by a family
of real polynomials $f_n\in \R[x]$, where
\[f_n(x)=r_{0,n}+r_{1,n}x+...+r_{m,n}x^m,\ \ \ n\in\N.\]

Indeed, what we really have  is a function  $f=[f_n]$, where\linebreak  $[f_n]:\R^*\mapsto \R^*$ is defined by
\begin{equation}\label{hp}[f_n]([(d_n)])=[(f_1(d_1),f_2(d_2),...)].\end{equation}

Set
\[I=\{n\in\N: f_n(a_n)<0\},  \ \  J=\{n\in\N: f_n(b_n)>0\}.\]

Since $I,J\in\mathcal U$, the intersection $I\cap J$ also belongs to $\mathcal U$. We can take into consideration only  real polynomials  $f_n$ with indexes $n\in I\cap J$. Thus,  for any $f_n$, where $n\in I\cap J$, we have
 \[f_n(a_n)\cdot f_n(b_n)<0.\]

 By Proposition 4, there exists $c_n\in (a_n,b_n)$ such that $f_n(c_n)=0$, for $n\in I\cap J$.
 For indices $n\in \N\setminus I\cap J$, we can take $c_n=0$.\footnote{The equality $[(r_n)]=[(s_n)]$ holds if $r_n=s_n$ for indices $n$, which belong  to some element  of the ultrafilter $\mathcal U$. Thus,  defining a hyperreal $[(c_n)]$ only those terms $c_n$ matter, of which indices $n$  belong to same element of the ultrafilter $\mathcal U$.}
 Put $c=[(c_n)]$.  By   (\ref{hp}), the equality $f(c)=0$ holds.

 Since
 \[a_n<c_n<b_n,\ \ \ \mbox{with}\ \ n\in I\cap J, \]
 we have the inequalities
 \[ [(a_n)]<[(c_n)]<[(b_n)].\]

 Hence $c\in (a,b)$.

  \hfill{$\Box$}

\pv To end the paper, by Propositions (1) and (5),  we obtain

\pv \textbf{Corollary.} The hyperreal number field  is real closed and the field  $\R^*(\sqrt{-1})$ is algebraically closed.

%\end{document}

\pvn{\tiny Piotr Błaszczyk,
Institute of Mathematics, Pedagogical University of Cracow, Poland\\
pb@up.krakow.pl}

\end{document}